\documentclass[a4paper,12pt]{amsart}

\usepackage[english]{babel}
\usepackage{amsthm}

\newcommand{\Set}[1]{\left\{\, #1 \,\right\}}
\newcommand{\Span}[1]{\langle\, #1 \,\rangle}
\newcommand{\Order}[1]{\lvert #1 \rvert}
\newcommand{\Size}[1]{\Order{#1}}

\DeclareMathOperator{\Sym}{Sym}
\DeclareMathOperator{\Alt}{Alt}
\DeclareMathOperator{\GF}{\mathbf{GF}}

\renewcommand{\phi}[0]{\varphi}
\renewcommand{\theta}[0]{\vartheta}
\renewcommand{\epsilon}[0]{\varepsilon}

\newtheorem{dummy}{Dummy}
\numberwithin{dummy}{section}
\numberwithin{equation}{section}

\newtheorem{theorem}[dummy]{Theorem}
\newtheorem{lemma}[dummy]{Lemma}

\newtheorem{cor}[dummy]{Corollary}

\theoremstyle{definition}

\theoremstyle{remark}

\newtheorem{fact}{Fact}

\begin{document}

\date{12 June 2006 --- Version 5.01%
}

\title[Imprimitivity of groups generated by round functions]
{Imprimitive  permutations groups\\
generated by the round  functions of\\
key-alternating block ciphers and \\
truncated differential cryptanalysis}

\author{A.~Caranti}

\address[A.~Caranti]{Dipartimento di Matematica\\
  Universit\`a degli Studi di Trento\\
  via Sommarive 14\\
  I-38050 Povo (Trento)\\
  Italy} 

\email{caranti@science.unitn.it} 

\urladdr{http://www-math.science.unitn.it/\~{ }caranti/}

\author{F.~Dalla Volta}

\address[F.~Dalla Volta]{Dipartimento di Matematica e Applicazioni\\
  Edificio U5\\
  Universit\`a degli Studi di Milano--Bicocca\\
  Via R.~Cozzi, 53\\
  I-20126 Milano\\
  Italy}

\email{francesca.dallavolta@unimib.it}

\urladdr{http://www.matapp.unimib.it/\~{ }dallavolta/}

\author{M.~Sala}

\address[M.~Sala]{Boole Centre for Research in Informatics\\
University College Cork\\
Cork\\
Ireland}

\email{msala@bcri.ucc.ie}

\author{Francesca Villani}

\begin{abstract}
 We answer a question of  Paterson, showing that all block systems for
 the group generated by the round functions of a key-alternating block
 cipher are the translates of a linear subspace.  Following up remarks
 of  Paterson  and  Shamir,  we  exhibit  a  connection  to  truncated
 differential cryptanalysis.

 We also give a condition  that guarantees that the group generated by
 the   round  functions   of   a  key-alternating   block  cipher   is
 primitive. This applies in particular to AES.
\end{abstract}


\thanks{Caranti   and  Dalla  Volta   are  members   of  INdAM-GNSAGA,
 Italy. Caranti  has been partially  supported by MIUR-Italy  via PRIN
 2001012275 ``Graded  Lie algebras  and pro-p-groups of  finite width,
 loop  algebras, and  derivations''.  Dalla  Volta has  been partially
 supported    by   MIUR-Italy    via   PRIN    ``Group    theory   and
 applications''.    Sala    has    been   partially    supported    by
 STMicroelectronics contract ``Complexity  issues in algebraic Coding
 Theory and Cryptography''}

\maketitle\thispagestyle{empty}

\section{Introduction}

Kenneth Paterson  \cite{Pat} has considered iterated  block ciphers in
which  the  group  generated  by the  one-round  functions  acts
imprimitively  on the  message space,  with the  aim of  exploring the
possibility  that this  might lead  to  the design  of trapdoors.  The
blocks  of imprimitivity  he uses  are  the translates  (cosets) of  a
linear subspace. He  asked whether it is possible  to construct other,
non-linear blocks of imprimitivity.

In the first part of this paper we answer this question in the
negative for key-alternating block ciphers, and exhibit a connection
to truncated differential cryptanalysis, following up remarks of
Paterson and Shamir.

We  then develop  a  conceptual  recipe to  guarantee  that the  group
generated by the one-round  functions of a key-alternating block
cipher  acts  primitively on  the  message  space.  We show  that  the
conditions we require are satisfied in a natural way by AES.

\subsection*{Acknowledgements}

We  are grateful to  P.~Fitzpatrick and  C.~Traverso for  their useful
comments. Part of this work has been presented at the Workshop on
Coding and Cryptography which was held at BCRI, UC Cork in 2005.

\section{Preliminaries}
\label{sec:prelim}

Let $G$ be a finite group, acting transitively on a set $V$. We write
the action of an element $g \in G$ on an element $\alpha \in V$
\emph{on the right}, that is, as
$\alpha g$. 
Also, $\alpha G = \Set{ \alpha g : g \in G}$ is the orbit
of $\alpha$ under $G$, and $G_{\alpha} = \Set {g \in G : \alpha g =
  \alpha }$ is the stabilizer of $\alpha$ in $G$.

A  \emph{partition}  of $V$  is  a  family  $\mathcal{B}$ of  nonempty
subsets of  $V$ such that  any element of  $V$ lies in  precisely one
element  of $\mathcal{B}$.  A partition  $\mathcal{B}$ is  said  to be
\emph{$G$-invariant} if for any $B  \in \mathcal{B}$ and $g \in
G$,  one  has  $B   g  \in  \mathcal{B}$.  A  $G$-invariant  partition
$\mathcal{B}$ is said  to be \emph{trivial} if $\mathcal{B}  = \Set{ V
}$, or $\mathcal{B} = \Set{ \Set{\alpha} : \alpha \in V}$.

A non-trivial, $G$-invariant partition of $V$ is said to be a \emph{block
system} for the action of $G$ on $V$. If such a block system exists,
then we say that $G$ is \emph{imprimitive} in its action on $V$
(equivalently, \emph{$G$ acts imprimitively} on $V$), \emph{primitive}
otherwise. An element $B$ of some block system $\mathcal{B}$ is called
a \emph{block}; since $G$ acts transitively on $V$, we have then $\mathcal{B}
= \Set{ B g : g \in G }$.

We note the following elementary
\begin{lemma}[\cite{Cam}, Theorem~1.7]
\label{lemma}
Let $G$ be a finite group, acting transitively on a set $V$.
Let $\alpha \in V$.

Then the blocks $B$ containing $\alpha$ are in one-to-one correspondence
with the subgroups $H$, with $G_{\alpha} < H < G$. The correspondence
is given by $B = \alpha H$.

In particular, $G$ is primitive if and only if $G_{\alpha}$ is a
maximal subgroup of $G$.
\end{lemma}

We will need  a fact from the basic theory of  finite fields. (See for
instance~\cite{Jac}~or \cite{LN}.)  Write $\GF(p^{n})$ for  the finite
field with $p^{n}$ elements, $p$ a prime.
\begin{lemma}
  \label{lemma:ff}
  $\GF(p^{n}) \subseteq \GF(p^{m})$ if and only if $n$ divides $m$.
\end{lemma}

In the rest of the paper, we tend to adopt the notation of~\cite{AES}.

Let $V =  V(n_{b}, 2)$, the vector space of dimension $n_{b}$ over the field
$\GF(2)$ with two elements, be the state space. $V$ has $2^{n_{b}}$ elements. 

For any $v \in V$, consider the translation by $v$, that is the map
\begin{equation*}
  \begin{aligned}
  \sigma_{v} : V &\to V,\\
  w &\mapsto w + v.
  \end{aligned}
\end{equation*} 
In particular, $\sigma_{0}$ is the identity map on $V$.  The set
\begin{equation*}
  T = \Set{\sigma_{v} : v \in V}
\end{equation*}
is an elementary abelian, regular subgroup of $\Sym(V)$. In fact, the map
\begin{equation}\label{eq:isoVT}
  \begin{aligned}
    V &\to T\\
    v &\mapsto \sigma_{v}\\
  \end{aligned}
\end{equation}
is an isomorphism of the additive group $V$ onto the multiplicative
group $T$.

We consider  a \emph{key-alternating block  cipher} (see Section~2.4.2
of~\cite{AES})  which consists of a
number  of  iterations   of  a  round  function  of   the  form  $\rho
\sigma_{k}$.   (Recall that  we  write maps  left-to-right, so  $\rho$
operates first.) Here  $\rho$ is a fixed permutation  operating on the
vector  space $V  =  V(n_{b},  2)$, and  $k  \in V$  is  a round
key. (According to the more general definition of~\cite{AES}, $\rho$
might depend on the round.)
Therefore each round consists of an application of $\rho$, followed by
a key addition.  This  covers for instance AES with \emph{independent}
subkeys.  Let $G  = \Span{\rho  \sigma_{k} :  k \in  V}$ the  group of
permutations of $V$ generated by the round functions. Choosing $k = 0$
we see  that $\rho \in G$,  and thus $T \le  G$. It follows  that $G =
\Span{T, \rho}$.

\section{Imprimitivity}
\label{sec:impr}


Kenneth Paterson  \cite{Pat} has considered iterated  block ciphers in
which  the  group  generated  by the  one-round  functions  acts
imprimitively  on the  message space,  with the  aim of  exploring the
possibility  that this  might lead  to  the design  of trapdoors.  The
blocks  of imprimitivity  he uses  are  the translates  (cosets) of  a
linear subspace. He  asked whether it is possible  to construct other,
non-linear blocks of imprimitivity:
\begin{quote}
  Can ``undetectable'' trapdoors based on more complex systems of
  imprimitivity be inserted in otherwise conventional ciphers?
  It is  easily shown that, in a  DES-like cipher, any [block]
  system  based  on a  linear  sub-space and  its  cosets  leads to  a
  noticeable regularity in the XOR  tables of small S-boxes.  It seems
  that we must look beyond  the ``linear'' systems considered here, or
  consider other types of round function.
\end{quote}
In a  personal communication  \cite{PatPers}, Paterson
remarks further
\begin{quote}
  At the  FSE conference  where it was  presented, Adi Shamir  told me
  that he could break the scheme using a truncated differential attack
  [\dots]
\end{quote}
Truncated differential cryptanalysis has been introduced in~\cite{Knu} by
L.~R.~Knudsen; see also the approach in~\cite{Wag}.

In this section we  answer Paterson's question for the key-alternating
block ciphers described above, by showing
\begin{theorem}
  \label{theo:impr}
  Let  $G$  be  the  group  generated  by the  round  functions  of  a
  key-alternating block cipher.  Suppose $G$ acts imprimitively on the
  message space.  Then the blocks of imprimitivity  are the translates
  of a linear subspace.
\end{theorem}

\begin{proof}
In the notation above, suppose $G$ acts imprimitively on $V$.

If $G$ has a nontrivial block  system, this is also a block system for
$T$. So if $\mathcal{B}$ is a block system for $G$, and $B \in \mathcal{B}$ is
the  block  containing $0$,  because of Lemma~\ref{lemma} we  have  $B
= 0  H$,  for some  $1 <  H  < 
T$. Because of the isomorphism \eqref{eq:isoVT}, we have
\begin{equation*}
  H = \Set{ \sigma_{u} : u \in U },
\end{equation*}
for a suitable subspace $U$ of $V$, with $U \ne \Set{0}, V$.
Since $T = \Set{\sigma_{v} : v \in V}$ is abelian, we have
\begin{multline*}
  \mathcal{B}
  =
  \Set{ B \sigma_{v} : v \in V }
  =
  \Set{ 0 H  \sigma_{v} : v \in V}
  =
  \Set{ 0 \sigma_{v} H : v \in V}
  =\\=
  \Set{ v H  : v \in V}
  =
  \Set{ v + U : v \in V}.
\end{multline*}
This completes  the proof  of the first  implication. The  converse is
immediate.
\end{proof}


\section{Truncated differential cryptanalysis}

We  now develop  a  relation to  truncated differential cryptanalysis,
elaborating on Shamir's comment.

Suppose $G$ acts  imprimitively on the message space  $V$, and use the
notation of the  proof of Theorem~\ref{theo:impr}. Let $v  \in V$. Now
$v H \rho$ is the block containing $v \cdot 1 \cdot \rho = v \rho$, so
that
\begin{equation*}
  v H \rho =  v \rho H,
\end{equation*}
for all
$v$. This  means that for  all $v  \in V$ and  $u \in U$  there is
$u' \in U$ such that
\begin{equation*}
  v \sigma_{u} \rho = (v + u) \rho 
  =
  v \rho + u' = v \rho \sigma_{u'}.
\end{equation*}
In other words we have the following connection to truncated
differential cryptanalysis.
\begin{cor}
  \label{fact:UinU}
  Suppose $G$ acts imprimitively on the message space $V$.

  Then there is a subspace $U \ne \Set{0}, V$ such that if $v, v + u
  \in V$  are two  messages whose difference  $u$ lies in  the subspace
  $U$, then the output difference also lies in $U$. 

  In other words, if
  $v \in V$ and $u \in U$, then
  \begin{equation}
    (v + u) \rho + v \rho \in U.
  \end{equation}

  Conversely, if the last condition holds, then $G$ acts imprimitively
  on $V$.
\end{cor}

To our  understanding, a subspace $U$  as in Corollary~\ref{fact:UinU}
could indeed be used as a  trapdoor as in Paterson's scheme, and still
be difficult to detect. This is  most clear when $U$ is chosen to have
dimension half of that of $V$.  To a cryptanalyst who knows $U$,
the complexity of  a brute force search is  reduced from $\Size{V}$ to
$2  \sqrt{\Size{V}}$.  However,  the number  of subspaces  of  a given
dimension $m$  of a finite vector  space of (even)  dimension $n$ over
$\GF(2)$ is largest for $m =  n/2$, and is $O(2^{m^{2}})$. If $U$
is not  just given by the vanishing  of some of the  defining bits, it
appears to us that it might be hard to find. Because of this, in the
next section  we approach the problem  of proving in  a conceptual way
that such a $U$ does not exists for a given key-iterated block cipher.

\section{Ensuring primitivity}

Ralph Wernsdorf has proved in~\cite{WAES} that the group $G$ generated
by  the  round  functions  of  AES with  independent  subkeys  is  the
alternating group $\Alt(n)$. Thus $G$ is definitely primitive on
$V$.  

In  the following  we  review this  consequence of  Wernsdorf's
result from a conceptual point of view. This comes in the form of
a  recipe  for  the  group  generated  by the  round  functions  of  a
key-alternating block cipher  to be primitive. We will  show that this
recipe is satisfied by AES in a rather natural way.

We begin with making the description of a key-alternating block cipher
we gave in Section~\ref{sec:prelim} more precise. (Again, we are
staying close to the notation of~\cite{AES}.) We assume $\rho = \gamma
\lambda$, where $\gamma$ and $\lambda$ are permutations. Here $\gamma$
is a bricklayer transformation, consisting of a number of S-boxes. The
message space $V$ is written as a direct sum
\begin{equation*}
  V = V_{1} \oplus \dots \oplus V_{n_{t}},
\end{equation*}
where each $V_{i}$ has the same dimension $m$ over $\GF(2)$. For $v
\in V$, we will write $v = v_{1} + \dots + v_{n_{t}}$, where $v_{i} \in
V_{i}$. Also, we consider the projections $\pi_{i} : V \to V_{i}$,
which map $v \mapsto v_{i}$. We have
\begin{equation*}
  v \gamma = v_{1} \gamma_{1} \oplus \dots \oplus v_{n_{t}} \gamma_{n_{t}},
\end{equation*}
where the $\gamma_{i}$ are S-boxes, which we allow to be different for
each $V_{i}$.

$\lambda$ is a linear mixing layer.

In AES  the S-boxes  are all  equal, and consist  of inversion  in the
field  $\GF(2^{8})$   with  $2^{8}$   elements  (see  later   in  this
paragraph), followed by an affine transformation.  The latter map thus
consists of  a linear transformation, followed by  a translation. When
interpreting AES  in our  scheme, we take  advantage of  the well-known
possibility of moving the linear  part of the affine transformation to
the linear  mixing layer, and  incorporating the translation in  the key
addition (see for instance~\cite{BES}). Thus  in our scheme for AES we
have $m = 8$, we identify  each $V_{i}$ with $\GF(2^{8})$, and we take
$x  \gamma_{i}  =  x^{2^{8}-2}$,  so that  $\gamma_{i}$  maps  nonzero
elements to their inverses, and zero to zero. As usual, we abuse
notation and write $x \gamma_{i} = x^{-1}$. Note, however, that
with this convention $x
x^{-1} = 1$ only for $x \ne 0$.

Our result, for a key-alternating block cipher as described earlier in
this section, is the following.
\begin{theorem}
  \label{theorem:main}
  Suppose the following hold:
  \begin{enumerate}
  \item  
    \label{item:01} 
    $0 \gamma = 0$ and $\gamma^{2} = 1$, the identity transformation.
  \item There is $1 \le r < m/2$ such that for all $i$
    \begin{itemize}
      \label{item:r}
    \item for all $0 \ne v \in V_{i}$, the image of the map $V_{i}
      \to V_{i}$,  which maps  $x \mapsto (x  + v) \gamma_{i}  + x
      \gamma_{i}$, has size greater than $2^{m-r-1}$, and
    \item there is no subspace of $V_{i}$, invariant
      under $\gamma_{i}$, of codimension less than or equal to $2 r$.
    \end{itemize}
  \item 
    \label{item:affine}
    No  sum of  some  of the  $V_{i}$  (except $\Set{0}$  and $V$)  is
    invariant under $\lambda$.
  \end{enumerate}
  Then $G$ is primitive.
\end{theorem}

We note immediately
\begin{lemma}
  \label{lemma:AES_is_OK}
  AES satisfies the hypotheses of Theorem~\ref{theorem:main}.
\end{lemma}

\begin{cor}
  \label{cor:AES_is_OK}
  The group generated by the round functions of AES with independent
  subkeys is primitive.
\end{cor}

\begin{proof}[Proof of Lemma~\ref{lemma:AES_is_OK}]
  Condition~\eqref{item:01}   is   clearly   satisfied.

  So  is~\eqref{item:affine},  by   the  construction  of  the  mixing
  layer. In fact,  suppose $U \ne \Set{0}$ is a  subspace of $V$ which
  is invariant  under $\lambda$. Suppose, without  loss of generality,
  that   $U   \supseteq   V_{1}$.   Because   of   \texttt{MixColumns}
  \cite[3.4.3]{AES}, $U$ contains the whole first column of the state.
  Now   the  action   of   \texttt{ShiftRows}  \cite[3.4.2]{AES}   and
  \texttt{MixColumns} on the first column shows that $U$ contains four
  whole columns, and  considering  (if the state has more than
  four  columns)  once  more  the  action  of  \texttt{ShiftRows}  and
  \texttt{MixColumns} one sees that $U = V$.

  The first part of  Condition~\eqref{item:r} is also well-known to be
  satisfied,  with  $r =  1$  (see~\cite{Ny}  but also~\cite{DR}).  We
  recall  the short  proof for  convenience. For  $a \ne  0$,  the map
  $\GF(2^{8}) \to  \GF(2^{8})$, which  maps $x \mapsto  (x +  a)^{-1} +
  x^{-1}$, has image  of size $2^{7} - 1$. In fact,  if $b \ne a^{-1}$,
  the equation
  \begin{equation}
    \label{eq:Ny1}
    (x + a)^{-1} + x^{-1} = b
  \end{equation}
  has at most two solutions. Clearly $x = 0, a$ are not solutions, so
  we can multiply by $x (x + a)$ obtaining the equation
  \begin{equation}
    \label{eq:Ny2}
    x^{2} + a x + a b^{-1} = 0,
  \end{equation}
  which has at most two solutions. If $b = a^{-1}$,
  equation~\eqref{eq:Ny1} has four solutions. Two of them are $x = 0,
  a$. Two more come from~\eqref{eq:Ny2}, which becomes
  \begin{equation*}
    x^{2} + a x + a^{2} 
    = 
    a^{2} \cdot\left( (x/a)^{2} + x/a + 1 \right)
    = 0.
  \end{equation*}
  By Lemma~\ref{lemma:ff}, $\GF(2^{8})$ contains $\GF(4) = \Set{0, 1,
  c, c^{2}}$, where $c, c^{2}$ are the roots of
  $y^{2} + y + 1 = 0$, Thus when $b = a^{-1}$ equation~\eqref{eq:Ny1}
  has the four solutions $0, a, a c, a c^{2}$. It follows that the
  image of the map $x \mapsto  (x +  a)^{-1} +
  x^{-1}$ has size
  \begin{equation*}
    \frac{2^{8} - 4}{2} + \frac{4}{4} = 2^{7} - 1,
  \end{equation*}
  as claimed.  

  As to  the second part  of Condition~\eqref{item:r}, one  could just
  use \textsf{GAP}~\cite{GAP4}  to  verify  that  the  only  nonzero  subspaces  of
  $\GF(2^{8})$ which are invariant under inversion are the subfields.
  According to  Lemma~\ref{lemma:ff}, the largest proper one is 
  thus $\GF(2^{4})$,  of codimension $4 > 2 = 2 r$.
  However,     this     follows      from     the     more     general
  Theorem~\ref{theorem:Sandro}, which we give in the Appendix.
\end{proof}

\begin{proof}[Proof of Theorem~\ref{theorem:main}]
  Suppose, by way of contradiction, that $G$ is imprimitive. According
  to Corollary~\ref{fact:UinU}, there is a subspace $U \ne \Set{0}, V$
  of  $V$ such  that  if $v,  v  + u  \in V$  are  two messages  whose
  difference $u$ lies in the  subspace $U$, then the output difference
  also lies in $U$, that is
  \begin{equation*}
    (v + u) \rho + v \rho \in U.
  \end{equation*}
  Since $\lambda$ is linear, we have
  \begin{fact}
    \label{fact:W}
    For all $u \in U$ and $v \in V$ we have
    \begin{equation}
      \label{eq:W}
      (v + u) \gamma + v \gamma \in U \lambda^{-1} = W,
    \end{equation}
  where $W$ is  also a linear subspace of $V$, with $\dim(W) = \dim(U)$.
  \end{fact}

  Setting $v = 0$ in~\eqref{eq:W}, and because of
  Condition~\eqref{item:01}, we obtain
  \begin{fact}
    \label{fact:UW}
    $U \gamma = W$ and $W \gamma = U$.
  \end{fact}

  Now if $U \ne \Set{0}$, we will have $U \pi_{i} \ne \Set{0}$ for
  some $i$. We prove some increasingly stronger facts under this
  hypothesis.

  \begin{fact}
    Suppose $U \pi_{i} \ne \Set{0}$ for some $i$. Then $W \cap V_{i}
    \ne \Set{0}$.
  \end{fact}
  Let  $u \in  U$,  with $u_{i}  \ne 0$.  Take  any $0  \ne v_{i}  \in
  V_{i}$. Then $(u + v_{i}) \gamma  + v_{i} \gamma \in W$, and also $u
  \gamma \in W$, by Fact~\ref{fact:UW}. It follows that $u \gamma + (u
  + v_{i})  \gamma + v_{i}  \gamma \in W$.  The latter vector  has all
  nonzero  components but  for the  one  in $V_{i}$,  which is  $u_{i}
  \gamma_{i} +  (u_{i} +  v_{i}) \gamma_{i} +  v_{i} \gamma_{i}  \in W
  \cap V_{i}$. If the latter vector is zero for all $v_{i} \in V_{i}$,
  then  the image  of the  map $V_{i}  \to V_{i}$,  which  maps $v_{i}
  \mapsto  (v_{i}   +  u_{i})  \gamma_{i}  +   v_{i}  \gamma_{i}$,  is
  $\Set{u_{i} \gamma_{i}}$,  of size  $1$. This contradicts  the first
  part of Condition~\eqref{item:r}.
  
  Clearly $(W \cap V_{i})\gamma = U \cap V_{i}$.  It follows
  \begin{fact}
    \label{fact:nonzerocap}
    Suppose $U \pi_{i} \ne \Set{0}$ for some $i$. Then $U \cap V_{i}
    \ne \Set{0}$.
  \end{fact}
  
  Finally we obtain
  \begin{fact}
    \label{fact:contains}
    Suppose $U \pi_{i} \ne \Set{0}$ for some $i$. Then $U \supseteq V_{i}$.
  \end{fact}
  According to Fact~\ref{fact:nonzerocap}, there is  $0 \ne u_{i} \in U
  \cap V_{i}$.  By the first part of  Condition~\eqref{item:r} the map
  $V_{i} \to V_{i}$, which maps $x  \mapsto (x + u_{i}) \gamma_{i} + x
  \gamma_{i}$, has  image of size  $> 2^{m-r-1}$. Since this  image is
  contained in the linear subspace $W \cap V_{i}$, it follows that the
  latter has size at least $2^{m-r}$, that is, codimension at most $r$
  in  $V_{i}$. The  same holds  for  $U \cap  V_{i} =  (W \cap  V_{i})
  \gamma$.   Thus  the linear  subspace  $U  \cap  W \cap  V_{i}$  has
  codimension at most $2 r$ in $V_{i}$. In particular, it is different
  from $\Set{0}$,  as $m  > 2 r$.  From Fact~\ref{fact:UW}  it follows
  that  $U \cap W  \cap V_{i}$  is invariant  under $\gamma$.   By the
  second part of Condition~\eqref{item:r} we have $U \cap W \cap V_{i}
  = V_{i}$, so that $U \supseteq V_{i}$ as claimed.
  
  From Fact~\ref{fact:contains} we obtain immediately
  \begin{fact}
    $U$ is a direct sum of some of the $V_{i}$, and $W = U$
  \end{fact}
  The second part follows from the fact that $W = U \gamma$, and
  $V_{i} \gamma = V_{i}$ for all $i$.
  
  Since $U =  W \lambda$ by~\eqref{eq:W}, we obtain $U  = U \lambda $,
  with     $U     \neq      \Set{0},     V$.      This     contradicts
  Condition~\eqref{item:affine}, and completes the proof.
\end{proof}

The proof of Theorem~\ref{theorem:main} can be adapted to prove a
slightly more general statement, in which Conditions~\eqref{item:01}
and \eqref{item:r} are replaced with
\begin{itemize}
  \item[($1'$)]  $0  \gamma =  0$  and $\gamma^{s}  = 1$,  for some $s > 1$.
  \item[($2')$] There is $1 \le r < m/s$ such that for all $i$
    \begin{itemize}
    \item[$\bullet$] for all $0 \ne v \in V_{i}$, the image of the map $V_{i}
      \to V_{i}$,  which maps  $x \mapsto (x  + v) \gamma_{i}  + x
      \gamma_{i}$, has size greater than $2^{m-r-1}$, and
    \item[$\bullet$] there is no proper subspace of $V_{i}$, invariant
	  under $\gamma_{i}$, of codimension less than or equal to $s r$.
    \end{itemize}
\end{itemize}

\section{Appendix}

We    are    grateful     to    Sandro    Mattarei (see~\cite{Matt},
and also~\cite{Others}, for more general results) for the following
\begin{theorem}
  \label{theorem:Sandro}
  Let $F$ be a field of characteristic two. Suppose  $U \ne 0$ is an
  additive subgroup of $F$ which contains the inverses of each of its
  nonzero elements. Then $U$ is a subfield of $F$.
\end{theorem}

\begin{proof}
  Hua's identity, valid in any associative (but not necessarily
  commutative) ring $A$, shows
  \begin{equation}
    \label{eq:Hua}
    a + ((a - b^{-1})^{-1} - a^{-1})^{-1} = a b a
  \end{equation}
  for $a, b \in A$, with $a, b, a b - 1$ invertible.

  First of all, $1 \in U$. This is because $U$ has even order, and
  each element different from $0,1$ is distinct from its inverse.

  Now~\eqref{eq:Hua} for $b = 1$, and $a \in U \setminus \Set{0, 1}$
  shows that for $a \in U$, also $a^{2} \in U$. (This is clearly valid
  also for $a = 0, 1$.) It follows that any $c \in U$ can be
  represented in the form $c = a^{2}$ for some $a \in
  U$. Now~\eqref{eq:Hua} shows that $U$ is closed under products, so
  that $U$ is a subring, and thus a subfield, of $F$.
\end{proof}

\providecommand{\bysame}{\leavevmode\hbox to3em{\hrulefill}\thinspace}
\providecommand{\MR}{\relax\ifhmode\unskip\space\fi MR }
\providecommand{\MRhref}[2]{%
  \href{http://www.ams.org/mathscinet-getitem?mr=#1}{#2}
}
\providecommand{\href}[2]{#2}

\end{document}